\documentclass [11pt]{article}
\begin{document}
\thispagestyle{empty}
\null\vspace{-1cm}
\medskip
\vspace{1.75cm}
\begin{center}

{\bf Signed  degree  sets  in  signed  graphs}\\
\vspace {1in}
{\bf S. Pirzada$^{*}$ , T. A. Naikoo$^{**}$  and  F. A. Dar$^{***}$}\\
\vspace{.8cm}
{\bf Department of Mathematics , University of Kashmir , Srinagar - 190006 , India }\\
\end{center}
{* Email :     sdpirzada@yahoo.co.in}\\
{** Email :     tariqnaikoo@rediffmail.com}\\
{*** Email:      sfarooqdar@yahoo.co.in}\\
\vspace {1in}
\begin{center}
{\bf Abstract}
\end{center}

The set D of distinct signed degrees of the vertices in a signed graph G is called its signed degree set . In this paper, we prove that every non-empty set of positive (negative) integers is the signed degree set of some connected signed graph and determine the smallest possible order for such a signed graph . We also prove that every non-empty set of integers is the signed degree set of some connected signed graph .

\bigskip
\vfill

\newpage

\vspace{0.15in}

\noindent{\bf 1. Introduction }\\
\vspace{0.15in}

All graphs in this paper are finite , undirected , without loops and multiple edges . A  signed graph G is a graph in which each edge is assigned a positive or a negative sign . These were first discovered by Harary [3] . The signed degree of a vertex $v_{i}$ in a signed graph G is denoted by  sdeg$(v_{i})$(or simply by $d_{i}$)and is defined as the number of positive edges incident with $v_{i}$less the number of negative edges incident with $v_{i}$.So, if$v_{i}$is incident with  $d_{i}^{+}$positive edges and $d_{i}^{-}$ negative edges, then sdeg$(v_{i})$= $d_{i}^{+}$-$d_{i}^{-}$.A signed degree sequence $\sigma=[d_{1},d_{2},...,d_{n}]$of a signed graph G is formed by listing the vertex signed degrees in non-increasing order . A sequence  $\sigma=[d_{1},d_{2},...,d_{n}$integers is graphical if  $\sigma$  is  a signed degree sequence of some signed graph . Also, a non-zero sequence $\sigma=[d_{1},d_{2},...,d_{n}$is a standard sequence if $\sigma$ is non-increasing ,$\sum_{i=1}^{n}d_{i}$ is even,$d_{1}>0$,each $|d_{i}|<n$, and $|d_{1}|\ge|d_{n}|$.\\ 

The following result , due to Chartrand et al. [1] , gives a necessary and sufficient condition for a sequence of integers to be graphical , which is similar to Hakimi's result for degree sequences [2] .\\ 
\vspace{0.3cm}

\noindent{\bf Theorem 1.1.}Let $\sigma=[d_{1},d_{2},...,d_{n}]$ be a standard sequence . Then , s  is graphical if and only if there exist integers  r  and  s  with $d_{1}=r-s$ and $0\le s\le\frac{n-1-d_{1}}{2}$such that\\
$$\sigma^{/}=[d_{2}-1,d_{3}-1,...,d_{r+1}-1,d_{r+2},d_{r+3},...,d_{n-s},d_{n-s+1}+1,...,d_{n}+1,]$$\\is graphical .\\

The next characterization for signed degrees in signed graphs is given by Yan et al. [5] . \\
\vspace{0.3cm}

\noindent{\bf Theorem 1.2.} A standard integral sequence $\sigma=[d_{1},d_{2},...,d_{n}$] is graphical if and only if \\
$$\sigma_{m}^{/}=[d_{2}-1,...,d_{d_{1}+m+1}-1,d_{d_{1}+m+2},...,d_{n-m},d_{n-m+1}+1,...,d_{n}+1,]$$\\
is graphical , where  m  is the maximum non-negative integer such that $d_{d_{1}+m+1}>d_{n-m+1}$\\

In [4], Kapoor et al. proved that every non-empty set of distinct positive integers is the degree set of a connected graph and determined the smallest order for such a graph .\\

\noindent{\bf 2. Main Results}
\vspace{0.15in}

First we have the following definition .\\
\noindent{\bf Definition.} The set D of distinct signed degrees of the vertices in a signed graph G is called its signed degree set .\\

	            Now, we obtain the following results .\\
\vspace{0.3cm}

\noindent{\bf Theorem 2.1.} Every non-empty set D of positive integers is the signed degree set of some connected signed graph and the minimum order of such a signed graph is N + 1, where N is the maximum integer in the set D.\\
{\bf Proof.} Let D be a signed degree set and $n_{0}$(D) denotes the minimum order of a signed graph G realizing D. Since N is the maximum integer in D, therefore there is a vertex  in G which is adjacent to at least N other vertices with a positive sign. Then,$n_{0}(D)\ge N + 1$. Now, if there exists a signed graph of order N + 1 with D as signed degree set , then $n_{0}(D)$ = N + 1. The existence of such a signed graph is obtained by using induction on the number of elements of D.\\

           Let $D = \{d_{1} , d_{2} ,~.~.~.~, d_{n}\},$ where $d_{1} < d_{2} <~.~.~.~< d_{n} ,$ be a set of positive  integers . For n = 1, let G be a complete graph on $d_{1} + 1$ vertices , that is $K_{d_{1}+1}$ in which each edge is assigned a positive sign. Then ,\\
$$sdeg(v) = (d_{1} + 1 - 1) - 0 = d_{1} , for~~all~~~ v\in V(G)$$\\
Therefore , G is a signed graph with signed degree set D = $\{d_{1}\}$.\\
	            
For n = 2 , let $G_{1}$ be a complete graph on $d_{1}$ vertices, that is $K_{d_{1}}$,in which each edge is assigned a positive sign and let $G_{2}$ be a null graph on $d_{2} - d_{1} + 1 > 0$ vertices , that is $K_{d_{2} - d_{1} + 1}$   Join every vertex of $G_{1}$ to each vertex of $G_{2}$ with a positive edge , so that we obtain a signed graph G on $d_{1} + d_{2} - d_{1} + 1 = d_{2} + 1$ vertices with\\
 $$sdeg(u) = (d_{1}- 1) + (d_{2} - d_{1} + 1) - 0 = d_{2} ,   for~~ all~~~ u\in V(G_{1}) ,$$\\
and\\
$$sdeg(v) = (0) + (d_{1}) - 0 = d_{1} ,   for~~ all~~~ v\in V(G_{2})$$.\\
Therefore , signed degree set of G is D = {$d_{1} , d_{2}$}.\\

For n = 3, let $G_{1}$ be a complete graph on $d_{1}$ vertices, that is $K_{d_{1}}$  , in which each edge is assigned a positive sign, $G_{2}$ be a complete graph on $d_{2} - d_{1} + 1 > 0$ vertices , that is $K_{d_{2} - d_{1} + 1 }$   in which each edge is assigned a positive sign, and $G_{3}$ be a null graph on $d_{3} - d_{2} > 0$ vertices, that is $\bar{K}_{d_{3} - d_{2}}$   Join every vertex of $G_{1}$ to each vertex of $G_{2}$ with a positive edge and join every vertex of $G_{1}$ to each vertex of $G_{3}$ with a positive edge, so that we obtain a signed graph G on $d_{1} + d_{2} - d_{1} + 1 + d_{3} - d_{2} = d_{3} + 1$ vertices with\\
$$sdeg(u) = (d_{1} - 1) + (d_{2} - d_{1} + 1) + (d_{3} - d_{2}) - 0 = d_{3} ,  for~~ all~~~ u\in V(G_{1}) ,$$\\   
$$sdeg(v) = (d_{2} - d_{1} + 1 - 1) + (d_{1}) - 0 = d_{2}  ,   for~~ all~~~ v\in V(G_{2}) ,$$\\
and\\
$$sdeg(w) = (0) + (d_{1}) - 0 = d_{1}  ,   for~~ all~~~ w\in V(G_{3}) .$$\\
Therefore , signed degree set of G is D = $\{d_{1} , d_{2} , d_{3}\}$.\\

Assume that the result holds for k . We show that the result is true for k + 1.\\

Let D = $\{d_{1} , d_{2},~.~.~.~, d_{k} , d_{k+1}\}$ be a k + 1 set of positive integers with $d_{1} < d_{2}  <~.~.~.~< d_{k}  < d_{k+1}$ . Clearly , $0 < d_{2} - d_{1} < d_{3} - d_{1} <~.~.~.~< d_{k} - d_{1}$ . Therefore , by induction hypothesis , there is a signed graph $G_{1}$ realizing the signed degree set $D_{1} = \{d_{2} - d_{1} , d_{3} - d_{1} ,~.~.~.~, d_{k} - d_{1}\}$ on $d_{k} - d_{1} + 1$ vertices as $|V(D_{1})| < k $. Let $G_{2}$ be a complete graph on $d_{1}$ vertices, that is $K_{d_{1}}$  , in which each edge is assigned a positive sign and $G_{3}$ be a null graph on $d_{k+1} - d_{k} > 0$ vertices, that is $\bar{K}_{d_{k+1} - d_{k}}$. Join every vertex of $G_{2}$ to each vertex of $G_{1}$ with a positive edge and join every vertex of $G_{2}$ to each vertex of $G_{3}$ with a positive edge, so that we obtain a signed graph G on $d_{k} - d_{1} + 1 + d_{1} + d_{k+1} - d_{k} = d_{k+1} + 1$ vertices with\\
$$sdeg(u) = (d_{i} - d_{1})+(d_{1}) - 0 = d_{i},for~~ all~~~ u\in V(G_{1})~~where~~~ 2\le i\le k ,$$\\    
$$sdeg(v) = (d_{1} - 1) + (d_{k} - d_{1}+1) + (d_{k+1} -  d_{k}) - 0 = d_{k+1} ,  for~~ all~~~ v\in V(G_{2}) ,$$\\
and\\
$$sdeg(w) = (0) + (d_{1}) - 0 = d_{1}  ,  for~~ all~~~ w\in V(G_{3}) .$$\\
Therefore , signed degree set of G is D = {$d_{1} , d_{2} ,~.~.~.~, d_{k} , d_{k+1}$} . Clearly , by construction , all the signed graphs are connected . Hence, the result follows .\\
\vspace{0.3cm}

\noindent{\bf Theorem 2.2.} Every non-empty set D of negative integers is the signed degree set          of some connected signed graph and the minimum order of such a graph is $|M|+ 1$ , where  M is the minimum integer in the set D .\\
{\bf Proof.} . Let D be a signed degree set and $m_{0}(D)$ denotes the minimum order of a signed graph G realizing D . Since $|M|$ is the maximum integer in D , therefore there is a   vertex in G which is adjacent to at least $|M|$ other vertices with a negative sign . Then ,$m_{0}(D)\ge|M|+ 1$. Now , if there exists a signed graph of order $|M|+ 1$ with D as signed degree set , then $m_{0}(D) = |M|+ 1 $.\\

           Let D = $\{-d_{1} , -d_{2} ,~.~.~.~, -d_{n}\} ,~~~ -d_{1} > -d_{2} >~.~.~.~> -d_{n}$ , be a set of negative integers where $d_{1} , d_{2} ,…, d_{n}$ are positive integers . Now , $D_{1} = \{d_{1} , d_{2} ,~.~.~.~, d_{n}\}$ is a set of positive integers with $d_{1} < d_{2} <….< d_{n}$ . By Theorem 2.1 , there exists a connected signed graph $G_{1}$ on $d_{n} + 1 = |-d_{n}| + 1$ vertices with signed degree set $D_{1} = \{d_{1} , d_{2} ,~.~.~.~, d_{n}\}$ . Now , construct a signed graph G from $G_{1}$ by interchanging positive edges with negative edges . Then , G is a connected signed graph on $|-d_{n}|+ 1$ vertices with degree set D = $\{-d_{1} , -d_{2} ,~.~.~.~, -d_{n}\}$. This proves the result .\\
\vspace{0.3cm}
                                            
\noindent{\bf Theorem 2.3.} Every non-empty set D of integers is the signed degree set of some connected signed graph .\\
{\bf Proof.}  Let D be a set of n integers . We have the following cases .\\
{\bf Case I.} D is a set of positive (negative) integers. Then , the result follows by Theorem 2.1( Theorem 2.2 ) .\\
{\bf Case II.} D = $\{0\}$. Then , a null graph G on one vertex, that is $K_{1}$, has signed degree set    D = $\{0\}$.\\
{\bf Case III.} D is a set of non-negative (non-positive) integers. Let $D = D_{1}\cup\{0\},$ where $D_{1}$ is a set of positive (negative) integers. Then , by Theorem 2.1 ( Theorem 2.2 ) , there is a signed graph $G_{1}$ with degree set $D_{1}$ . Let $G_{2}$ be a null graph on two vertices , that is $\bar{K}_{2}$  . Let e = uv be an edge in $G_{1}$ with positive (negative) sign and let $x , y\in V(G_{2})$ . Add the positive (negative) edges ux and vy , and the negative (positive) edges uy and vx , so that we obtain a connected signed graph G with signed degree set D . We note that addition of such edges do not effect the signed degrees of the vertices of $G_{1}$ , and the vertices x and y have signed degrees zero each .\\

{\bf Case IV.} D is a set of non-zero integers . Let $D = D_{1}\cup D_{2}$ , where $D_{1}$ is a set of positive integers and $D_{2}$ is a set of negative integers . Then , by Theorem 2.1 and Theorem 2.2 , there are connected signed graphs $G_{1}$ and $G_{2}$ with signed degree sets $D_{1}$ and $D_{2}$ . Let $e_{1} = uv$ be an edge in $G_{1}$ with positive sign and $e_{2} = xy$ be an edge in $G_{2}$ with negative sign . Add the positive edges ux and vy , and the negative edges uy and vx , so that we obtain a connected signed graph G with signed degree set D . We note that addition of such edges do not effect the signed degrees of the vertices of $G_{1}$ and $G_{2}$ .\\

{\bf Case V.} D is a set of integers . Let $D = D_{1}\cup D_{2}\cup\{0\}$ , where $D_{1}$ and $D_{2}$ are the sets of positive and negative integers respectively . Then, by Theorem 2.1 and Theorem 2.2 , there are connected signed graphs $G_{1}$ and $G_{2}$ with signed degree sets $D_{1}$ and $D_{2}$ . Let $G_{3}$ be a null graph on one vertex , that is $K_{1}$ . Let $e_{1} = uv$ be an edge in $G_{1}$ with positive sign , and let $x\in V(G_{2})$ and $y\in V(G_{3})$ . Add the positive edges uy and vx , and the negative edges ux and vy , so that we obtain a connected signed graph G with signed degree set D. We note that addition of such edges do not effect the signed degrees of the vertices of $G_{1}$ and $G_{2}$ , and the vertex y has signed degree zero . This completes the proof .\\

\vspace{0.3cm}

\noindent{\bf Theorem 2.4.}. If G is a signed graph with vertex set V where $|V|= r$ and signed degree set $\{d_{1} , d_{2} ,~.~.~.~, d_{n}\}$ . Then , for each $k\ge 1$, there is a signed graph with kr vertices and signed degree set $\{d_{1} , d_{2} ,~.~.~.~, d_{n}\}$ .\\                   

{\bf Proof.} For each i , $1\le i\le k$ , let $G_{i}$ be a copy of G with vertex set $V_{i}$ . Define a signed graph H with vertex set $W =\cup_{i=1}^{k}v_{i}$   where $V_{i}\cap V_{j} =\phi(i\neq j)$ and the edges of H are the edges of $G_{i}$ for all i, where $1\le i\le k$ . Therefore , H is a signed graph on kr vertices with signed degree set $\{d_{1} , d_{2} ,~.~.~.~, d_{n}\}$.\\

\vspace{0.15in}

\noindent{\bf REFERENCES}

\vspace{0.15in}

[1]  G. Chartrand , H. Gavlas , F. Harary and M. Schultz , On signed degrees in signed 

       graphs , Czeck. Math. J. 44 (1994) 677-690.(Zbl   0837.05110)\\

[2]  S. L. Hakimi , On the realizability of a set of integers as degrees of the vertices of 

       a graph , SIAM  J. Appl. Math. 10 (1962) 496-506.(Zbl   0109.16501)\\

[3]  F. Harary, On the notion of balance in a signed graph , Michigan Math. J. 2(1953) 

       143-146.(Zbl   0056.42103)\\

[4]  S. F. Kapoor, A. O. Polimeni and C. E. Wall, Degree sets for graphs, Fund. Math. 

       65 (1977) 189-194.(Zbl   0351.05129)\\

[5]   J. H. Yan , K. W. Lih , D. Kuo  and  G. J. Chang , Signed  degree  sequences  of    

       signed graphs , J. Graph Theory 26 (1997) 111-117.(Zbl   0890.05072)

\end{document}